\documentstyle[12pt]{article}
\begin{document}
{\bf Journal's headings:}
\vskip 7cm
\begin{center}
{\large \bf SMALL PARAMETER PERTURBATIONS OF NONLINEAR PERIODIC SYSTEMS
$^\dagger$}
\end{center}
\begin{center}
Mikhail Kamenski$^1$, Oleg Makarenkov$^1$ and Paolo Nistri$^2$
\end{center}
\begin{center}
$^1$Department of Mathematics,\\
 Voronezh State University,
Voronezh, Russia.\\
E-mail: Mikhail\verb+@+kam.vsu.ru,$\;\;$ E-mail: !omakarenkov\verb+@+kma.vsu.ru!\\
$^2$ Dipartimento di Ingegneria dell'Informazione,\\
Universit\`a di Siena, 53100 Siena, Italy.\\
E-mail: pnistri\verb+@+dii.unisi.it
\end{center}

\date{}
\def\thefootnote{\fnsymbol{footnote}}
\footnotetext[2]{Supported by the national research project MIUR:
``Feedback control and optimal control''!,
by RFBR! grants 02-01-00189 and 02-01-00307 !and U.S.CRDF - RF Ministry of Education grant VZ-010!}
\def\thefootnote{\fnsymbol{footnote}}

\def\inte{\int\limits}     \let\qq=\qquad       \let\q=\quad
\let\w=\widetilde          \let\wh=\widehat     \let\mx=\mbox
\let\ol=\overline          \let\D=\Delta        \let\d=\delta
\let\e=\epsilon            \let\g=\gamma        \def\mm{{-\!\circ}}
\let\G=\Gamma              \let\ba=\beta        \def\dist{{\fam0 dist\,}}
\let\a=\alpha              \let\th=\theta       \let\nn=\nonumber
\let\s=\sigma              \let\O=\Omega      \def\L{{\cal L}}
\def\N{\mbox{\bf N}}       \def\bp{\mbox{\bf P}} \let\emp=\emptyset
\def\Z{\mbox{\bf Z}}       \def\ind#1{\mathop{#1}\limits}
\def\C{\mbox{\bf C}}       \let\l=\lambda       \let\sm=\setminus
\let\o=\omega              \def\meas{{\fam0 meas\,}}
\def\nor#1{{\left\|\,#1\,\right\|}}             \let\bcup=\bigcup
\def\R{\mbox{\bf R}}       \let\vf=\varphi    \def\sca#1#2{\langle #1,#2\rangle}
\def\Re{{\rm Re\,}}         \def\Im{{\rm Im\,}}        \def\DD{{\cal D}}
\let\ds=\displaystyle      \let\lra=\longrightarrow
\def\dfrac#1#2{\ds{#1\over #2}} \let\y=\eta
\let\sb=\subset            \def\J{{\cal J}}
\def\qed{\hfill $\squar$}
\def\squar{\vbox{\hrule\hbox{\vrule height 6pt \hskip
6pt\vrule}\hrule}}
\def\toto{{\begin{array}{c}\rightarrow\\*[-2ex]
             \rightarrow \end{array}}}
\catcode`\@=11
\def\section{\@startsection {section}{1}{\z@}{-3.5ex plus-1ex minus
-.2ex}{2.3ex plus.2ex}{\bf}}

\noindent {\bf Abstract.}$\;\;$ In this paper we consider a class of nonlinear
periodic differential
systems perturbed by two nonlinear periodic terms with multiplicative different
powers of a small
parameter $\varepsilon>0$. For such a class of systems we provide conditions
which guarantee the existence
of periodic solutions of given period $T>0$. These conditions are expressed in
terms of the behaviour
on the boundary of an open bounded set $U$ of ${\rm R}^n$ of the solutions of
suitably defined
linearized systems. The approach is based on the classical theory of the
topological degree for
compact vector fields. An application to the existence of periodic solutions to
the van der Pol
equation is also presented.



\vspace{6mm}
\noindent
\section {Introduction}

\vspace{2mm}
\noindent
This paper represents the continuation of the work started by the authors in
\cite{kmn}.
There we considered the existence problem both of periodic solutions and of
solutions
of the Cauchy problem for a system of ordinary differential equations described
by
\begin{equation} \label{1'}
  \dot x=\psi(t,x)+\varepsilon\phi(t,x),
\end{equation}
where $\phi,\ \psi:{\rm R}\times{\rm R}^n\to {\rm R}^n$
are continuously differentiable, $T$-periodic with
respect to time $t$, functions and $\varepsilon$ is a small positive parameter.

\noindent
To solve this problem for $\varepsilon>0$ sufficiently small a new approach
was presented in \cite{kmn}. Such an approach is based on the linearized system
\begin{equation} \label{2'}
\dot y=\frac{\partial\psi}{\partial x}(t,\Omega(t,0,\xi))y +
\phi(t,\Omega(t,0,\xi)),
\end{equation}
where $\xi\in {\rm R}^n$ and $\Omega(\cdot,t_0,\xi)$ denotes the solution of
(\ref{1'}) at
$\varepsilon=0$ satisfying $x(t_0)=\xi$. Specifically, consider the change of
variable
 \begin{equation} \label{3'}
  z(t)= \Omega(0,t,x(t)),
\end{equation}
and the solution $\eta(\cdot,s,\xi)$ of (\ref{2'}) such that $y(s)=0$. If there
exists a bounded
open set $U\subset {\rm R}^n$ such that $\Omega(T,0,\xi)=\xi$ for any
$\xi\in\partial U$, and
$\eta(T,s,\xi)-\eta(0,s,\xi)\not=0,$ for any $s\in[0,T],$  and any $\xi\in
\partial U.$
Then (\ref{1'}) has a $T$-periodic solution for $\varepsilon>0$ sufficiently
small provided
that $\gamma(\eta(T,0,\cdot))\not=0$. Here $\gamma(F,U)$ denotes the rotation
number of a continuous
map $F:\overline U\to \overline U$.

\noindent
The advantage of the proposed approach as compared with the classical averaging
method, which is one of the most
useful tool for treating this problem, mainly consists
in the fact that in order to use this second method for establishing the
existence of periodic solutions
in perturbed systems of the form (\ref{1'}) one must assume that the change of
variable (\ref{3'})
is $T$-periodic with respect to $t$ for every $T$-periodic function $x$ such
that $\Omega(0,t, x(t))\in U,
t\in [0,T]$, instead that only on the boundary of the bounded open set $U$.

\noindent
The same assumption is necessary in vibrational control problems, \cite{s} and
\cite{b}, to reduce the considered
system to the standard form for applying the averaging method.
For an extensive list of references on this topic see \cite{bullo}.

\noindent
Our approach has been also employed in \cite{kmn1}
to prove the existence of periodic solutions for a class of first order
singularly perturbed differential systems.

\noindent
The aim of this paper is to extend the previously outlined approach to a more
general class of
perturbed systems than (\ref{1'}).
Precisely, we consider here systems of the following form:
\begin{equation} \label{1bis}
  \dot
x=\psi(t,x)+\varepsilon^2\phi_1(t,x)+\varepsilon^3\phi_2(t,x,\varepsilon),
\end{equation}
where the functions $\psi,\ \phi_1:{\rm R}\times{\rm R}^n\to {\rm R}^n,\
\phi_2:{\rm R}\times{\rm R}^n\times[0,1]\to {\rm R}^n$ are !continuously
differentiable! and $T$-periodic
with respect to the first variable
and $\varepsilon$ is a small positive parameter.
We denote again by $\Omega(\cdot,t_0,\xi)$ the solution of the Cauchy problem
\begin{eqnarray} \label{zam}
  \left\{\begin{array}{l}
            \dot x=\psi(t,x), \\
            x(t_0)=\xi,
         \end{array} \right.
\end{eqnarray}
and by $\eta_i(\cdot,s,\xi), i=1,2,$ the solution of the Cauchy problems
\begin{eqnarray} \label{aux}
  \left\{\begin{array}{l}
            \dot y=\dfrac{\partial \psi}{\partial x}(t,\Omega(t,0,\xi))y+
            \phi_1(t,\Omega(t,0,\xi)), \quad {\rm if}\ i=1,\\\\
            \dot y=\dfrac{\partial \psi}{\partial x}(t,\Omega(t,0,\xi))y+
            \phi_2(t,\Omega(t,0,\xi),0), \quad {\rm if}\ i=2, \\\\
            y(s)=0.
         \end{array} \right.
\end{eqnarray}

\noindent
In Section 2 we prove the main result of the paper:
!Theorem 1. Indeed, under
suitable assumptions on $\Omega(T,0,\xi)$ and $\eta_i(T,s,\xi), i=1,2,$ for
$s\in [0,T]$
and $\xi \in \partial U$, for $\varepsilon >0$ sufficiently small, we prove the
existence of $T$-periodic solutions of system (\ref{1bis}) provided that
${\rm deg}(\eta_2(T,0,\cdot),U)\not= 0$. Moreover, as it is shown in Theorem 2,
in the
case of system (\ref{1'}) this result implies Theorem 1 of \cite{kmn}.!

\vskip.2truecm
\noindent
In Section 3, we apply Theorem 1 to autonomous systems in ${\rm R}^2$ perturbed
by a
non-autonomous term of higher order (with respect to $\varepsilon>0$) to show
the existence of periodic solutions in a suitable open set defined by means of
the
trajectory of the linear part of the autonomous system.

\noindent
Finally, in Section 4 we illustrate the obtained results by means of an example
concerning the van der Pol equation. To the best knowledge of the authors this
represents
a new approach to investigate the existence of periodic solutions for the
periodically forced
van der Pol equation. Indeed, many papers in the literature are devoted to the
study of the response
of the van der Pol equation to periodic stimulus (of different period),
but the methods are quite different, in fact they are essentially based on
asymptotic
expansions, Fourier series, singular !perturbation! theory and averaging
methods. We refer to the papers
\cite{belyakova}, \cite{braaksma}, \cite{el-abbasy}, \cite{grammel} and the
references therein.

\section{Main results}

Consider the system
\begin{equation} \label{1}
  \dot
x=\psi(t,x)+\varepsilon^2\phi_1(t,x)+\varepsilon^3\phi_2(t,x,\varepsilon),
\end{equation}
where the functions $\psi,\ \phi_1:{\rm R}\times{\rm R}^n\to {\rm R}^n,\
\phi_2:{\rm R}\times{\rm R}^n\times[0,1]\to {\rm R}^n$
are continuously differentiable, $T$-periodic with
respect to time $t$, and $\varepsilon$ is a small positive parameter.

\vskip0.2truecm
\noindent
To investigate the existence of $T$-periodic solutions of system (\ref{1}) we
introduce
the compact integral operator
$F_{\varepsilon}: C([0,T], {\rm R}^n) \to  C([0,T], {\rm R}^n)$ defined by
$$
F_{\varepsilon}(x)(t)=x(t)-x(T)-\int\limits_0^t\left(\psi(\tau,x(\tau))+\varepsilon^2\phi_1(\tau,x(\tau))+
\varepsilon^3\phi_2(\tau,x(\tau),\varepsilon)\right)d\tau,
$$
for any $t\in[0,T]$. Clearly if $F_{\varepsilon}(x)=0$ then $x$ is a
$T$-periodic solution of system (\ref{1}).

\vskip0.2truecm
\noindent
Denote by $\Omega(\cdot,t_0,\xi)$ the solution of the Cauchy problem
\begin{eqnarray} \label{zam}
  \left\{\begin{array}{l}
            \dot x=\psi(t,x), \\
            x(t_0)=\xi,
         \end{array} \right.
\end{eqnarray}
and by $\eta_i(\cdot,s,\xi), 1=1,2,$ the solution of the Cauchy problems
\begin{eqnarray} \label{aux}
  \left\{\begin{array}{l}
            \dot y=\dfrac{\partial \psi}{\partial x}(t,\Omega(t,0,\xi))y+
            \phi_1(t,\Omega(t,0,\xi)), \quad {\rm if}\ i=1, \\\\
            \dot y=\dfrac{\partial \psi}{\partial x}(t,\Omega(t,0,\xi))y+
            \phi_2(t,\Omega(t,0,\xi),0), \quad {\rm if}\ i=2, \\\\
            y(s)=0.
         \end{array} \right.
\end{eqnarray}
\vskip0.2truecm
\noindent
The following lemma provides an explicit representation of the functions
$\eta_1$ and $\eta_2.$

\vskip0.2truecm
\noindent
{\bf Lemma 1.} !!Let $\xi\in{\rm R}^n$ and $s,t\in[0,T]$. We have that!!
\begin{eqnarray}
\eta_1(t,s,\xi)&=&\frac{\partial\Omega}{\partial
z}(t,0,\xi)\int\limits_s^t\Phi_1(\tau,\xi)d\tau \nonumber
\end{eqnarray}
where
\begin{eqnarray}
& & \Phi_1(t,\xi)=\frac{\partial\Omega}{\partial
z}(0,t,\Omega(t,0,\xi))\phi_1(t,\Omega(t,0,\xi)) \nonumber
\end{eqnarray}
and
\begin{eqnarray}
  \eta_2(t,s,\xi)&=&\frac{\partial\Omega}{\partial
z}(t,0,\xi)\int\limits_s^t\Phi_2(\tau,\xi,0)d\tau \nonumber
\end{eqnarray}
where
\begin{eqnarray}
& & \Phi_2(t,\xi,\varepsilon)=\frac{\partial\Omega}{\partial
z}(0,t,\Omega(t,0,\xi))
\phi_2(t,\Omega(t,0,\xi),\varepsilon).\nonumber
\end{eqnarray}

\noindent
{\bf Proof.}

\noindent
It is sufficient to observe, that the matrix $\dfrac{\partial \Omega}{\partial
z}(t,0,\xi)$ is the
fundamental matrix of the linear system
$$
\dot y=\frac{\partial \psi}{\partial x}(t,\Omega(t,0,\xi))y
$$
Moreover, $\left(\dfrac{\partial \Omega}{\partial
z}(t,0,\xi)\right)^{-1}=\dfrac{\partial\Omega}{\partial z}
(0,t ,\Omega(t,0,\xi)).$ In fact, if we derive with respect to $\xi$ the
identity
$$
\Omega(0,t, \Omega(t,0, \xi))=\xi
$$
we obtain
$$
\frac{\partial \Omega}{\partial z}(0,t,\Omega(t,0,\xi)) \frac{\partial
\Omega}{\partial z}
(t,0,\xi) = I,
$$
whenever $\xi\in{\rm R}^n$.
Therefore, by the variation of constants formula for linear nonhomogeneous
system
$$
 \dot y=\frac{\partial \psi}{\partial
x}(t,\Omega(t,0,\xi))y+\phi_1(t,\Omega(t,0,\xi)),
$$
we have
\begin{eqnarray}
  \eta_1(t,s,\xi)&=&\int\limits_s^t\frac{\partial\Omega}{\partial z}(t,0,\xi)
  \left(\frac{\partial \Omega}{\partial
z}(\tau,0,\xi)\right)^{-1}\phi_1(\tau,\Omega(\tau,0,\xi))d\tau=\nonumber\\
  &=&\frac{\partial\Omega}{\partial
z}(t,0,\xi)\int\limits_s^t\Phi_1(\tau,\xi)d\tau. \nonumber
\end{eqnarray}

\noindent
The formula for $\eta_2$ is obtained in the same way. \qed

\vskip0.4truecm
\noindent
We have the following result.

\vskip0.2truecm
\noindent
{\bf Theorem 1.}
Let $U\subset {\rm R}^n$ be an open and bounded set.
Assume, that

\vskip0.1truecm
(A1) $\Omega(T,0,\xi)=\xi,\ $ !for any! $\xi\in\partial U,$

(A2) $\eta_1(T,s,\xi)-\eta_1(0,s,\xi)=0,$ !for any! $s\in[0,T]$ !and any!
$\xi\in \partial U,$

(A3) $\eta_2(T,s,\xi)-\eta_2(0,s,\xi)\not=0,$ !for any! $s\in[0,T]$ !and any!
$\xi\in \partial U.$

\vskip0.2truecm
\noindent
Then for $\varepsilon>0$ sufficiently small
\begin{equation}\label{estab}
  {\rm deg}(F_{\varepsilon},W(T,U))={\rm deg}(\eta_2(T,0,\cdot),U),
\end{equation}
where $W(T,U)=\left\{x\in C\left([0,T],{\rm R}^n\right):\
\Omega(0,t,x(t))\in U,\ \mbox{whenever} \ t\in[0,T]\right\}.$

\vskip0.5truecm
\noindent
In order to prove the Theorem we need the following Lemma.

\vskip0.2truecm
\noindent
{\bf Lemma 2.}
Let $z\in C^1([0,T], {\rm R}^n), \; f\in C([0,T], {\rm R}^n)$ and $b\in {\rm
R}^n$. If

\begin{equation}\label{6}
\int_0^t \frac{\partial \Omega}{\partial z}(s,0,z(s)) \dot z(s) ds + z(0)=
b + \int_0^t f(s)ds,
\end{equation}
then
\begin{equation}\label{7}
z(t)= b + \int_0^t \frac{\partial \Omega}{\partial z}(0,s, \Omega(s,0,
z(s)))f(s) ds.
\end{equation}

\noindent
{\bf Proof.}

\noindent
Take the derivative of (\ref{6}) with respect to $t$ and then apply
$\dfrac{\partial \Omega}{\partial z}(0,t,\Omega(t,0,z(t)))$ to both sides.
Finally, integrating the resulting differential system from $0$ to $t$ and
observing that from
(\ref{6}) we have $z(0)=b$ one has (\ref{7}).\qed

\vskip0.4truecm
\noindent
{\bf Proof of Theorem 1.}

\noindent
Let $x$ be a solution of the equation
\begin{equation}\label{F0}
  F_{\varepsilon}(x)=0.
\end{equation}
Thus $x$ is a $T$-periodic solution to (\ref{1}). Consider the change of
variable
\begin{equation}\label{zam}
x(t)=\Omega(t,0,z(t)), \qquad t\in[0,T],
\end{equation}
with inverse given by
\begin{equation}\label{zam1}
z(t)=\Omega(0,t,x(t)), \qquad t\in[0,T].
\end{equation}
Observe that if $x$ is a solution of (\ref{F0}) then it is differentiable; in
fact, for any $t\in[0,T]$, we have
\begin{equation}\label{ei}
  x(t)=x(T)+\int_0^t\left(\psi(\tau,x(\tau))+\varepsilon^2\phi_1(\tau,x(\tau))+
\varepsilon^3\phi_2(\tau,x(\tau),\varepsilon)\right)d\tau.
\end{equation}
Therefore, from (\ref{zam1})
$z$ is also differentiable. Consider
\begin{equation}\label{16}
\frac{d}{dt} \Omega(t,0,z(t)) = \frac{\partial \Omega}{\partial t}(t,0,z(t))+
\frac{\partial \Omega}{\partial z}(t,0,z(t))\dot z(t),
\end{equation}
since
$$
\frac{\partial \Omega}{\partial t}(t,0,z(t))=\psi(t,\Omega(t,0,z(t)))
$$
from (\ref{16}) we have that
$$
\Omega(t,0,z(t)) - z(0) = \int_0^t \psi(s,\Omega(s,0,z(s)))ds +
\int_0^t \frac{\partial \Omega}{\partial z}(s,0,z(s))\dot z(s)ds,
$$
or equivalently,
\begin{equation}\label{17}
\Omega(t,0,z(t)) -  \int_0^t \psi(s,\Omega(s,0,z(s)))ds =
z(0) +  \int_0^t \frac{\partial \Omega}{\partial z}(s,0,z(s))\dot z(s)ds.
\end{equation}
By using (\ref{zam}), (\ref{17}) and Lemma 2 with
$b=\Omega(T,0,z(T))$ we can rewrite (\ref{ei}) in the following
form
\begin{equation}\label{G0}
  G_{\varepsilon}\left(z\right)(t)=0,
\end{equation}
where $G_{\varepsilon}: C([0,T], {\rm R}^n) \to  C([0,T], {\rm
R}^n)$ is given by
\begin{eqnarray}
G_\varepsilon (z)(t)=z(t)-\Omega(T,0,z(T))-\int
\limits_0^t \big(\varepsilon^2 \Phi_1(\tau,z(\tau))+
\varepsilon^3\Phi_2(\tau,z(\tau),\varepsilon)\big)d\tau. \nonumber
\end{eqnarray}
Therefore the solutions of the equation (\ref{F0}) belonging to the set $W(T,U)$
correspond to the solutions of the equation (\ref{G0}) belonging to the set
$$
Z=\{z \in C([0,T],{\rm R}^n):\ z(t) \in U,\ \mbox{whenever} \ t \in [0,T]\}.
$$
and by the homeomorphism Theorem for compact vector fields
(see \cite{kz}, Theorem 26.4) to prove Theorem 1 it is enough to
show that
\begin{equation}\label{indG}
{\rm deg}(G_{\varepsilon},Z)={\rm deg}(\eta_2(T,0,\cdot),U)
\end{equation}
for $\varepsilon>0$ sufficiently small.

\noindent For this, consider the compact vector field $\tilde
G_{\varepsilon}:  C([0,T],{\rm R}^n) \to C([0,T],{\rm R}^n)$
defined as follows
$$
 \tilde G_{\varepsilon}=I-A_{\varepsilon},
$$
where
$$
  A_{\varepsilon}(z)(t)=\Omega(T,0,z(T))+
    \int\limits_0^T\left(\varepsilon^2\Phi_1(\tau,z(\tau))+
    \varepsilon^3\Phi_2(\tau,z(\tau),\varepsilon)\right) d\tau,
$$
for any $t\in[0,T]$, !hence $A_{\varepsilon}(z)$ is a constant function
in $C([0,T],{\rm R}^n)$.!
Let us show that for $\varepsilon>0$
sufficiently small the compact vector fields $G_\varepsilon$ and
$\tilde G_{\varepsilon}$ are homotopic on the boundary of the set
$Z.$ To this aim, we define, for $\lambda \in [0,1]$ the following homotopy:
$$
\Delta_{\varepsilon}(\lambda,z)(t)=z(t)-\Omega(T,0,z(T))-
\int \limits_0^{\lambda t+(1-\lambda)T}\left(\varepsilon^2\Phi_1(\tau,z(\tau))+
\varepsilon^3\Phi_2(\tau,z(\tau),\varepsilon)\right) d\tau,
$$
whenever $t\in[0,T]$, which deforms the vector field
$G_\varepsilon$ to the vector field $\tilde G_{\varepsilon}$. Let
us show that $\Delta_{\varepsilon}$ does not vanish on the
boundary of the set $Z$ for $\varepsilon>0$ sufficiently small.

\noindent Assume the contrary. Therefore there exists a sequence
${\left\{\varepsilon_k\right\}}_{k=1}^{\infty}\subset(0,1]$ such
that $\varepsilon_k\to 0$ as $k\to\infty$ and sequences
${\left\{\lambda_k\right\}}_{k=1}^{\infty}\subset[0,1]$ and
${\left\{z_k\right\}}_{k=1}^{\infty}\subset \partial Z$ such that
\begin{eqnarray}
 z_k(t) & = & \Omega(T,0,z_k(T))+\varepsilon_k^2\int \limits_0^{\lambda_k t+
 (1-\lambda_k)T}\Phi_1(\tau,z_k(\tau))d\tau+ \nonumber \\
& &  +\varepsilon_k^3\int \limits_0^{\lambda_k t+
(1-\lambda_k)T}\Phi_2(\tau,z_k(\tau),\varepsilon_k)d\tau,\ \
t\in[0,T]. \label{f15}
\end{eqnarray}
Without loss of generality we can assume that
$\lambda_k\to\lambda_0$ and $z_k\to z_0\ $  in $C([0,T], {\rm
R}^n)$ as $k\to\infty.$ Therefore $\lambda_0\in[0,1]$ and
$z_0\in\partial Z.$ Furthermore, there exists a sequence
${\left\{t_k\right\}}_{k=1}^{\infty}$ such that
$z_k(t_k)\in\partial U$ and by condition (A1)
\begin{equation}\label{f21}
\Omega(T,0,z_k(t_k))=z_k(t_k),\ !\mbox{for any}!\ \ k\in{\rm N}.
\end{equation}
By subtracting from (\ref {f15}), where $t$ is replaced by $T$, the
same equation with $t=t_k$ we obtain
\begin{eqnarray}\label{f22}
  z_k(T)-z_k(t_k) & = & \varepsilon_k^2\int \limits_{\lambda_k t_k+
  (1-\lambda_k)T}^T\Phi_1(\tau,z_k(\tau))d\tau+ \nonumber\\
 & & + \varepsilon_k^3\int \limits_{\lambda_k t_k+
 (1-\lambda_k)T}^T\Phi_2(\tau,z_k(\tau),\varepsilon_k)d\tau.
\end{eqnarray}
By (\ref{f21}) the equation (\ref{f15}) where $t=T$ can be
rewritten as
\begin{eqnarray}
 & & z_k(T)-z_k(t_k)  = \Omega(T,0,z_k(T))-\Omega(T,0,z_k(t_k))+\nonumber\\
 & &  +\varepsilon_k^2\int \limits_0^T\Phi_1(\tau,z_k(\tau))d\tau+
 \varepsilon_k^3\int \limits_0^T\Phi_2(\tau,z_k(\tau),\varepsilon_k)d\tau.
\nonumber
\end{eqnarray}
or equivalently
\begin{eqnarray}\label{f23}
& &  \left(I-\frac{\partial \Omega}{\partial
z}(T,0,z_k(t_k))\right)\left(z_k(T)-z_k(t_k)\right)= \nonumber\\
& &  =\frac{\partial^2 \Omega}{\partial
z^2}(T,0,z_k(t_k))(z_k(T)-z_k(t_k))(z_k(T)-z_k(t_k)) + \nonumber\\
& &  +\varepsilon_k^2\int \limits_0^T\Phi_1(\tau,z_k(\tau))d\tau+
\varepsilon_k^3\int
\limits_0^T\Phi_2(\tau,z_k(\tau),\varepsilon_k)d\tau+\nonumber\\
& &  +o(z_k(t_k),z_k(T)-z_k(t_k)),
\end{eqnarray}
where the function $o(\xi,h)$ satisfy
\begin{equation}\label{o}
  \frac{\|o(\xi,h)\|}{\|h\|^2} \to 0,\ {\rm as}\ \|h\| \to 0,\quad \mbox{with}
\quad h,\ \xi\in {\rm R}^n.
\end{equation}
Replacing (\ref {f22}) in (\ref {f23}) and dividing by
$\varepsilon_k^3>0$ after a suitable transformation we obtain
\begin{equation}\label{weob}
 \frac{1}{\varepsilon_k}Q_{1k}(z_k)+Q_{2k}(z_k)=P_k+\frac{1}{\varepsilon_k^3}\
o(z_k(t_k),z_k(T)-z_k(t_k))
\end{equation}
where\\
\vskip0.01cm
\hskip1cm
$P_k=\varepsilon_k \dfrac{\partial^2 \Omega}{\partial z^2}(T,0,z_k(t_k))\cdot$\\

\vskip0.01cm
\hskip1cm $ \cdot \left(\int \limits_0^{\lambda_k
t+(1-\lambda_k)T}\Phi_1(\tau,z_k(\tau))d\tau+
\varepsilon_k\int \limits_0^{\lambda_k
t+(1-\lambda_k)T}\Phi_2(\tau,z_k(\tau),\varepsilon_k)d\tau\right)\cdot$ \\
\vskip0.01cm
\hskip1cm$ \cdot \left(\int \limits_0^{\lambda_k
t+(1-\lambda_k)T}\Phi_1(\tau,z_k(\tau))d\tau+
\varepsilon_k\int \limits_0^{\lambda_k
t+(1-\lambda_k)T}\Phi_2(\tau,z_k(\tau),\varepsilon)d\tau\right),$\\
and
\begin{eqnarray}
  Q_{1k}(z_k) & = & \left(I-\frac{\partial\Omega}{\partial
z}(T,0,z_k(t_k))\right)
  \int \limits_{\lambda_k
t_k+(1-\lambda_k)T}^T\Phi_1(\tau,z_k(\tau))d\tau-\nonumber\\
  & & -\int \limits_0^T\Phi_1(\tau,z_k(\tau))d\tau,\nonumber
\end{eqnarray}
\begin{eqnarray}
  Q_{2k}(z_k) & = & \left(I-\frac{\partial\Omega}{\partial
z}(T,0,z_k(t_k))\right)
  \int \limits_{\lambda_k
t_k+(1-\lambda_k)T}^T\Phi_2(\tau,z_k(\tau),\varepsilon_k)d\tau-\nonumber\\
 & &  -\int \limits_0^T\Phi_2(\tau,z_k(\tau),\varepsilon_k)d\tau.\nonumber
\end{eqnarray}

\noindent
It is easy to see, that
\begin{equation}\label{P}
  P_k\to 0\ {\rm as}\ k\to\infty
\end{equation}
and
\begin{equation}\label{o0}
  \frac{1}{\varepsilon_k^3}\, o(z_k(t_k),z_k(T)-z_k(t_k))\to 0\ {\rm as}\
k\to\infty.
\end{equation}

\noindent
Let us show, that
\begin{equation}\label{Q}
  \frac{1}{\varepsilon_k}Q_{1k}(z_k)\to 0\ {\rm as}\ k\to\infty.
\end{equation}

\noindent
We have
\begin{eqnarray} \label{Q1}
 \frac{1}{\varepsilon_k}Q_{1k}(z_k) & = &
\frac{Q_{1k}(z_k)-Q_{1k}(c_k)+Q_{1k}(c_k)}{\varepsilon_k}=\nonumber\\
& = & \frac{\partial Q_{1k}}{\partial
z}(c_k)\left(\frac{z_k-c_k}{\varepsilon_k}\right)+
\frac{Q_{1k}(c_k)}{\varepsilon_k}+\frac{o_k(z_k-c_k)}{\varepsilon_k},
\end{eqnarray}
where by $c_k$ we denote the constant function $c_k(t)\equiv
z_k(t_k)$ for any $t\in[0,T]$ and the function $o_k(\cdot)$
satisfy
\begin{equation}\label{o}
  \frac{\|o_k(h)\|}{\|h\|} \to 0,\ {\rm as}\ \|h\| \to 0,\ h\in {\rm
C}([0,T],{\rm R}^n).
\end{equation}

\noindent
By using (\ref{f22})
\begin{equation}
  \|z_k(t)-c_k\|\le \varepsilon_k^2 M,
\end{equation}
for any $t\in[0,T]$, where $M>0$ is a constant and so the first
and the third term in (\ref{Q1}) tends to zero as $k$ tends to
$\infty.$ Let us prove, that
\begin{equation}\label{toprove}
  Q_{1k}(c_k)=0,\ k\in{\rm N}.
\end{equation}
By Lemma 1
\begin{equation}\label{enough}
  \left(I-\frac{\partial\Omega}{\partial z}(T,0,\xi)\right)
  \int \limits_s^T\Phi_1(\tau,\xi)d\tau-\int \limits_0^T\Phi_1(\tau,\xi)d\tau=
  \eta_1(0,s,\xi)-\eta_1(T,s,\xi)
\end{equation}
and by condition (A2) we obtain (\ref{toprove}). Thus (\ref{Q})
holds true and by (\ref{P}), (\ref{o0}) and (\ref{Q}) we can pass
to the limit in (\ref{weob}), obtaining
\begin{equation}\label{Q20}\nonumber
  Q_2(z_0)=0,
\end{equation}
where $Q_2=\lim_{k\to\infty}Q_{2k}.$ On the other hand by Lemma 1
we obtain a result the analogous to (\ref{enough}) for $\Phi_2$, thus
\begin{equation}\label{enough2}\nonumber
  Q_2(z_0)=\eta_2(0,s,z_0(t_0))-\eta_2(T,s,z_0(t_0)),
\end{equation}
where $s=\lim_{k\to\infty}(\lambda_k t_k+(1-\lambda_k)T)$, hence
\begin{equation}\label{enough2}\nonumber
  \eta_2(0,s,z_0(t_0))-\eta_2(T,s,z_0(t_0))=0,
\end{equation}
which contradicts assumption (A3) since $z_0(t_0)\in\partial U.$
Therefore, there exists $\varepsilon_0>0$ such that
\begin{equation}\label{defor}
  \Delta_{\varepsilon}(\lambda,z)\not=0,\ \lambda\in[0,1],\ z\in\partial Z,\
\varepsilon\in(0,\varepsilon_0)
\end{equation}
and so
\begin{equation}\label{Gdeg}
  {\rm deg}(G_{\varepsilon},Z)={\rm deg}(\tilde G_{\varepsilon},Z),\
\varepsilon\in(0,\varepsilon_0).
\end{equation}

\noindent
Denote by $C_0([0,T],{\rm R}^n)$ the subspace of the space
$C([0,T],{\rm R}^n)$ consisting of all constant functions
defined on $[0,T]$ with values in ${\rm R}^n$.
We have $A_{\varepsilon}(\partial Z)\subset C_0([0,T],{\rm R}^n).$
From (\ref{defor}) with $\lambda=0$ we have
$$
  \tilde G_{\varepsilon}(z)\not= 0,\  z \in \partial Z, \
\varepsilon\in(0,\varepsilon_0),
$$
and using the inclusion
$\partial (Z\bigcap C_0([0,T],{\rm R}^n))\subset\partial Z$
we conclude that the vector field $\tilde G_{\varepsilon}$ does not have zeros
on the
boundary of the set $Z\bigcap C_0([0,T],{\rm R}^n)$ when
$\varepsilon\in(0,\varepsilon_0).$
Therefore by the reduction property of the topological degree, (see \cite{kz},
Theorem 27.1),
we have
\begin{equation} \label{vr1}
  {\rm deg}_{C([0,T],{\rm R}^n)}(\tilde G_{\varepsilon},Z)=
    {\rm deg}_{C_0([0,T],{\rm R}^n)}(\tilde G_{\varepsilon},
       Z\bigcap C_0([0,T],{\rm R}^n)).
\end{equation}
Furthermore, since the constant function $z\in Z\bigcap C_0([0,T],{\rm R}^n)$
is a solution of the equation $\tilde G_{\varepsilon}(z)=0$ if and only if
the element $\xi\in U,\ $ !given by!  $\xi=z(t),\ $ !for any! $ \ t\in [0,T],$
is a solution
of the equation $K_{\varepsilon}\xi=0,$ !where $K_{\varepsilon}: U \to {\rm
R}^n$ is defined
as follows!
$$
  K_{\varepsilon}\xi=\xi-\Omega(T,0,\xi)-\varepsilon^2\int
\limits_0^{T}\Phi_1(\tau,\xi)d\tau-\varepsilon^3\int
\limits_0^{T}\Phi_2(\tau,\xi,\varepsilon)d\tau,
$$
then
\begin{equation}\label{vr2}
  {\rm deg}_{C_0([0,T],{\rm R}^n)}(\tilde G_{\varepsilon},Z\bigcap
C_0([0,T],{\rm R}^n))=
    {\rm deg}_{{\rm R}^n}(K_{\varepsilon},U).
\end{equation}
Consider now !$K_{\varepsilon, 2}: U \to {\rm R}^n$ defined by!
$$
K_{\varepsilon,2} \xi=-\varepsilon^3\int
\limits_0^{T}\Phi_2(\tau,\xi,\varepsilon)d\tau.
$$
By condition (A1), i.e. $\xi-\Omega(T,0,\xi)=0,\ $ for any $ \ \xi\in\partial
U$,
Lemma 1 and assumption (A2) we obtain
\begin{eqnarray}
  -\int\limits_0^T\Phi_1(\tau,\xi)d\tau&=&+\int\limits_T^0\Phi_1(\tau,\xi)d\tau
    =\eta_1(0,T,\xi)=\nonumber\\
  &=&\eta_1(0,T,\xi)-\eta_1(T,T,\xi)=0,\ !\mbox{whenever}! \ \xi\in\partial
U.\nonumber
\end{eqnarray}
Therefore
$K_{\varepsilon}\xi=K_{\varepsilon,2} \xi,\ $ !whenever! $ \ \xi\in\partial U$
and thus
\begin{equation}\label{vr3}
  {\rm deg}_{{\rm R}^n}(K_{\varepsilon},U)=
    {\rm deg}_{{\rm R}^n}(K_{\varepsilon,2},U),\ !\mbox{for any}!
\varepsilon\in(0,\varepsilon_0).
\end{equation}
Finally, consider the vector field !$\tilde K_{\varepsilon, 2}: U \to {\rm R}^n$
defined as follows!

$$
\tilde K_{\varepsilon, 2} \xi=-\int_0^T\Phi_2(\tau,\xi,\varepsilon)d\tau.
$$
By condition (A2) we have
\begin{equation}\label{bc}
  -
\int\limits_0^T\Phi_2(\tau,\xi,0)d\tau=\eta_2(0,T,\xi)-\eta_2(T,T,\xi)\not=0,\
!\mbox{for any}! \ \xi\in\partial U.
\end{equation}
Therefore there exists $\varepsilon_1\in(0,\varepsilon_0)$ such that
$$
  - \int\limits_0^T\Phi_2(\tau,\xi,\varepsilon)d\tau\not=0,\ !\mbox{for any}! \
\xi\in\partial U \
!\mbox{and any}! \ \varepsilon\in(0,\varepsilon_1).
$$
and so the vector fields $K_{\varepsilon,2}$ and $\tilde K_{\varepsilon,2}$ are
linearly
homotopic on the boundary of the set $U$ for any
$\varepsilon\in(0,\varepsilon_1).$
Moreover, by condition (\ref{bc}) the vector fields $\tilde K_{\varepsilon,2}$
and $\tilde K_{0, 2}$
are linearly homotopic on $\partial U$ for $\varepsilon\in(0,\varepsilon_1).$
Therefore
\begin{eqnarray}
 & & {\rm deg} (K_{\varepsilon,2},U)={\rm
deg}(\eta_2(0,T,\xi)-\eta_2(T,T,\xi),U)=\nonumber\\
 & & = {\rm deg}(\eta_2(T,0,\xi)-\eta_2(0,0,\xi),U)={\rm
deg}(\eta_2(T,0,\xi),U).\label{Q1e}
\end{eqnarray}

\noindent
Summarizing (\ref{Gdeg}), (\ref{vr1}), (\ref{vr2}), (\ref{vr3}) and (\ref{Q1e})
we obtain
(\ref{estab}), !which is the claim of Theorem 1.! \qed

\vskip0.4truecm
\noindent
Consider now the system
\begin{equation} \label{1_}
  \dot x=\psi(t,x)+\varepsilon\phi(t,x,\varepsilon),
\end{equation}
and
$$
F_{\varepsilon}(x)(t)=x(t)-x(T)-\int\limits_0^t\left(\psi(\tau,x(\tau))+\varepsilon\phi(\tau,x(\tau),\varepsilon)\right)d\tau.
$$
Denote by $\eta$ the solution of Cauchy problem (\ref{aux}) associated with
(\ref{1_}).
As a direct consequence of the previous result we have the following.

\vskip0.2truecm
\noindent
{\bf Theorem 2.}
Let $U\subset {\rm R}^n$ be an open and bounded set.
Assume, that

\vskip0.1truecm
(A1) $\Omega(T,0,\xi)=\xi,\ $ !for any! $ \ \xi\in\partial U,$

(A2) $\eta(T,s,\xi)-\eta(0,s,\xi)\not=0,$ !for any! $ \ s\in[0,T],$ \ !and any!
$ \ \xi\in \partial U.$

\vskip0.2truecm
\noindent
Then for $\varepsilon>0$ sufficiently small
\begin{equation}\label{estab_}
  {\rm deg}(F_{\varepsilon},W(T,U))={\rm deg}(\eta(T,0,\cdot),U).
\end{equation}

\vskip0.4truecm
\noindent
Observe that Theorem 1 of \cite{kmn} follows from Theorem 2 above.

\section{Small perturbations of autonomous systems}

\noindent
In this Section by using the previous results we provide sufficient conditions
to ensure that the topological degree of the integral operator associated to the
following
system in ${\rm R}^2$
\begin{equation}\label{P2}
\dot x=Ax+\varepsilon^2\phi_1(x)+\varepsilon^3\phi_2(t,x,\varepsilon),
\end{equation}
is different from zero for $\varepsilon>0$ sufficiently small.

\noindent
System (\ref{P2}) is regarded here as the perturbed system of the autonomous
system
\begin{equation}\label{np}
\dot x=Ax+\varepsilon^2\phi_1(x)
\end{equation}
by means of the $T$-periodic perturbation
$\varepsilon^3\phi_2(t,x,\varepsilon)$,
where $A$ is a $2\times 2$ matrix. At $\varepsilon=0$ we have the linear system
\begin{equation}\label{ls}
\dot x=Ax.
\end{equation}
We assume, that

\vskip0.2truecm
\noindent
(A1)$\;$ !the matrix $A$ has eigenvalues $i \lambda $ and $-i \lambda $, where
$\lambda>0.$

\vskip0.1truecm
\noindent
Let $x_0$ be a !non-zero! periodic solution of system (\ref{ls}) which is the boundary of
an open
!! set $U_0$ of ${\rm R}^2$. Observe that such a periodic solution exists in
virtue
of (A1). Moreover, since the period $T$ of any periodic solution of system
(\ref{ls}) is equal to
$\frac{2\pi}{\lambda},$ we assume that the function $\phi_2$
is $\frac{2\pi}{\lambda}$-periodic with respect to the first variable.

\noindent
Denote by $U_{\delta}$ the open !! sets whose boundaries are given by the
trajectories
$(1+\delta)x_0$ for $\delta\in{\rm R}.$! Since $x_0$ is a
$\frac{2\pi}{\lambda}$-periodic solution of
system (\ref{ls}) then $(1+\delta)x_0$, $\delta\in{\rm R},$ is also a
$\frac{2\pi}{\lambda}$-periodic solution of this system.
Moreover we have $U_0\subset U_\delta.$

\noindent
Define $\Omega,\ $ $\eta_1,\ $ $\eta_2$ and $F_{\varepsilon}$ as in the previous
Section
with $\psi(t,\xi)=A\xi$ and $T=\frac{2\pi}{\lambda}.$

\vskip0.4truecm
\noindent
We have the following result.
\vskip0.2truecm
\noindent
{\bf Theorem 3.}
Assume condition (A1). Moreover, assume that for $\delta\in(0,1)$ we have

(A2) $\eta_1(T,s,\xi)-\eta_1(0,s,\xi)=0, \ $ !for any! $ \ s\in[0,T]\ $ !and
any! $\ \xi\in \partial U_0,$

(A3) $\eta_2(T,s,\xi)-\eta_2(0,s,\xi)\not=0, \ $ !for any!  $ \ s\in[0,T]\ $
!and any!  $\ \xi\in \partial U_0,$

(A4) $\eta_1(T,s,\xi)-\eta_1(0,s,\xi)\not=0, \ $ !for any! $ \ s\in [0,T]\ $
!and any! $\ \xi\in \partial U_{\delta}.$

\vskip0.2truecm
\noindent
Then for all $\varepsilon>0$ sufficiently small
\begin{eqnarray}
& & {\rm deg}\left(F_{\varepsilon},W\left(T,U_\delta\right)\backslash \overline
W\left(T,U_0\right) \right)
=\nonumber\\
& & =\left({\rm deg}\left(\eta_1\left(T,0,\cdot\right),U_{\delta}\right)-{\rm
deg}\left(\eta_2\left(T,0,\cdot\right),U_0\right)\right).\label{degF1}
\end{eqnarray}

\vskip0.2truecm
\noindent
{\bf Proof.}

\noindent
First of all observe that the set $W\left(T,U_\delta\right)
\backslash \overline W\left(T,{U}_0\right) $
is well defined since $U_\delta \supset U_0$.

\noindent
Moreover, from (A1), (A2) and (A3) it follows that Theorem 1 is applicable with
$U=U_0$ and so for $\varepsilon>0$ sufficiently small we obtain
\begin{equation}\label{ind1}\nonumber
  {\rm deg}(F_{\varepsilon},W(T,U_0))={\rm
deg}\left(\eta_2\left(T,0,\cdot\right),U_0\right).
\end{equation}
Finally, (A1) and (A4) imply that Theorem 2 is applicable with
$U=U_{\delta}$ and so for $\varepsilon>0$ sufficiently small we obtain
\begin{equation}\label{ind2}\nonumber
  {\rm deg}(F_{\varepsilon},W(T,U_\delta))={\rm
deg}\left(\eta_1\left(T,0,\cdot\right),U_{\delta}\right).
\end{equation}
Thus (\ref{ind1}) and (\ref{ind2}) ensure (\ref{degF1}) with
$\varepsilon>0$ sufficiently small. \qed

\vskip0.5truecm
\noindent
!We will show now that if we perturb (\ref{np}) by a $T_\varepsilon$-
periodic
perturbation, where 
\begin{equation}\label{near}
  T_\varepsilon\to T {\ \rm as\ } \varepsilon\to 0, 
\end{equation}
then (\ref{P2}) has
$T_\varepsilon$-periodic
solutions.
In case all the powers of $\varepsilon$ in the right hand side of (\ref{P2})
are the same and $T_\varepsilon\not=T$ this result is called 
the phenomenon of frequency pulling (see \cite{an}).
So we would like to show that in case of the perturbation (\ref{P2}),
where the powers of $\varepsilon$ are different,
this phenomenon 
 still has a place.
As $T_\varepsilon$ we consider here 
$T_{\varepsilon,\mu}=\frac{2\pi}{\lambda(1+\varepsilon^3\mu)},$
where $\mu$ is a scalling parameter.
More precisely we consider the system!
\begin{equation}\label{R2_}
  \dot
x=Ax+\varepsilon^2\phi_1(x)+\varepsilon^3\phi_2\left(\frac{t}{1+\varepsilon^3\mu
},x,\varepsilon,\mu\right).
\end{equation}
!Denote by $F_{\varepsilon, \mu}$ the integral operator associated with
(\ref{R2_}) and 
denote by $\eta_i(\cdot,s,\xi),\ $ $i=1,2$,  the solution of (\ref{aux})
corresponding
to (\ref{R2_}) for $\mu=0$.!

\vskip0.4truecm
\noindent
We can formulate the following result.

\vskip0.2truecm
\noindent
{\bf Theorem 4.}
Assume that for $\mu=0$ and  !$T=\frac{2\pi}{\lambda}$! all the conditions of Theorem 3 are
satisfied. Then there exists
$\mu_0>0$ such that for every $\mu\in(-\mu_0,\mu_0)$ there is an
$\varepsilon_0>0$ such that
\begin{eqnarray}
& & {\rm deg}\left(F_{\varepsilon, \mu},W\left(!T_{\varepsilon,\mu}!,
U_\delta\right)\backslash
\overline W\left(!T_{\varepsilon,\mu}!,U_0\right) \right) =\nonumber\\
& & =\left({\rm
deg}\left(\eta_1\left(!T_{\varepsilon,\mu}!,0,\cdot\right),U_{\delta}\right)-{\rm
deg}\left(\eta_2\left(!T_{\varepsilon,\mu}!,0,\cdot\right),U_0\right)\right) 
\label{degF4}
\end{eqnarray}
for any $\varepsilon\in(0,\varepsilon_0).$

\vskip0.2truecm
\noindent
{\bf Proof.}

\noindent
In (\ref{R2_}) consider the change of variable:
\begin{equation}\label{zam4}
  y(t)=x\left(t(1+\varepsilon^3\mu)\right).
\end{equation}
Thus we obtain the system
\begin{eqnarray}\label{aux1}
   \dot y &=& Ay+\varepsilon^2\phi_1(y)+\varepsilon^3\phi_{2,
\mu}(t,y,\varepsilon),
\end{eqnarray}
where
$$
  \phi_{2, \mu}(t,\xi,\varepsilon)=\phi_2(t,\xi,\varepsilon,\mu)+\mu
A\xi+\varepsilon^2\mu \ \phi_1(\xi)+\varepsilon^3\mu\
\phi_2(t,\xi,\varepsilon,\mu).
$$
Since the conditions of Theorem 3 hold true for (\ref{aux1}) for $\mu=0$
then there exists $\mu_0>0$ such that the
conditions of Theorem 3 hold true for (\ref{aux1}) for every fixed
$\mu\in(-\mu_0,\mu_0).$
Observe, that the solution of problem (\ref{aux}) with $i=1$
associated with (\ref{R2_}) and that one associated with (\ref{aux1}) coincide.
Denote by $\eta_{2, \mu}$ the solution of the problem (\ref{aux}) with $i=2$
associated
with (\ref{aux1}) and denote by $G_{\varepsilon, \mu}$ the compact operator
corresponding to the $T_0$-periodic problem for the system (\ref{aux1}).
By Theorem 3 for every $\mu\in(-\mu_0,\mu_0)$ there exists $\varepsilon_0>0$
such that
\begin{eqnarray}
& & {\rm deg}\left(G_{\varepsilon, \mu},W\left(T_0,U_\delta\right)\backslash
\overline W\left(T_0,U_0\right) \right) =\nonumber\\
& & =\left({\rm deg}\left(\eta_1\left(T_0,0,\cdot\right),U_{\delta}\right)-{\rm
deg}\left(\eta_{2, \mu}\left(T_0,0,\cdot\right),U_0\right)\right) \label{degF4_}
\end{eqnarray}
for $\varepsilon\in(0,\varepsilon_0).$
Since
$$
  \eta_{2, \mu}(T_0,0,\xi)\to\eta_2(T_0,0,\xi)\ {\rm as}\ \mu\to 0,\ !\mbox{for
any}! \ \xi\in{\rm R}^2
$$
we may assume without loss of generality that $\mu_0>0$ is choosen in such a way
that
\begin{equation}\label{wl}
  {\rm deg}\left(\eta_{2, \mu}\left(T_0,0,\cdot\right),U_0\right)=
{\rm deg}\left(\eta_2\left(T_0,0,\cdot\right),U_0\right),\ !\mbox{for any}! \
\mu\in(-\mu_0,\mu_0).
\end{equation}
By (\ref{zam4}) the zeros of the operator $G_{\varepsilon, \mu}$ belonging to
$W(T_0,U)$
correspond to the zeros of the operator $F_{\varepsilon, \mu}$
belonging to $W(!T_{\varepsilon,\mu}!,U)$. Therefore
\begin{eqnarray}
& & {\rm deg}\left(F_{\varepsilon, \mu},W\left(!T_{\varepsilon,\mu}!,U_\delta\right)
\backslash \overline W\left(!T_{\varepsilon,\mu}!, U_0\right) \right) =\nonumber\\
& & = {\rm deg}\left(G_{\varepsilon, \mu},W\left(T_0,U_\delta\right)
\backslash \overline W\left(T_0, U_0\right) \right)\nonumber
\end{eqnarray}
and by taking into account (\ref{degF4_}) and (\ref{wl}) we obtain
(\ref{degF4}).

\qed

\section{An application to the existence of periodic solutions of the van der
Pol equation}

We end the paper by illustrating a topological degree approach, based on the
above results, to
investigate the existence of periodic solutions of the van der Pol equation
\begin{equation}\label{vdp}
  \ddot x-\varepsilon(1-x^2)\dot x+x=0.
\end{equation}

\noindent
Consider the compact operator $F_\varepsilon: {\rm C}\left([0,2\pi],{\rm
R}^2\right) \to
{\rm C}\left([0,2\pi],{\rm R}^2\right)$ given by
\begin{eqnarray}\label{cov}
  F_\varepsilon(x)(t)=x(t)-x(2\pi)-\int\limits_0^t\left(
  \begin{array}{l}
   \qquad \quad x_2(\tau)\\
 - x_1(\tau) +\varepsilon(1-x_1^2(\tau))x_2(\tau)
   \end{array} \right) d\tau,
\end{eqnarray}
whose zeros $x=(x_1,x_2)$ correspond to $2\pi$-periodic
solution $(x,\dot x)$ of system (\ref{vdp}).
It is known (see for instance \cite{mi}, \S\ 4.6), that for $\varepsilon>0$
sufficiently small the operator
$F_\varepsilon$ has a zero $x_\varepsilon \in C([0,2\pi],{\rm R}^2)$ such that
$$ \|x_\varepsilon\|_{C_2\pi}\to 2\ {\rm as}\ \varepsilon\to 0. $$
Here $\|\cdot\|_{C_2\pi}$ denotes the norm in the Banach space $C([0,2\pi],{\rm
R}^2)$.
\vskip0.1truecm
\noindent
We have the following result.

\vskip0.2truecm
\noindent
{\bf Proposition 1.}
For $\varepsilon>0$ sufficiently small and $\delta\in(0,2)$
\begin{equation}
  {\rm ind}(F_\varepsilon,W(!2\pi!,U_\delta))=0,
\end{equation}
where
$$
  U_\delta=\left\{\xi\in{\rm R}^2:\ \|\xi\|\in(2-\delta,2+\delta)\right\}.
$$

\vskip0.2truecm
\noindent
{\bf Proof.}

\noindent
The proposition is a straightforward consequence of Theorem 2.
In fact, if we put
$\psi(t,\xi)=\left(\begin{array}{l}
\; \xi_2 \\  - \xi_1
\end{array} \right)$
and
$\phi(t,\xi,\varepsilon)=\left(
  \begin{array}{l}
    \qquad 0\\
(1-\xi_1^2)\xi_2
\end{array} \right)$
then, as it can be verified, we have that
$$
 \Omega(t,t_0,\xi)=\left(
  \begin{array}{cc}
    \cos (t-t_0) & \sin (t-t_0) \\ - \sin (t-t_0) & \cos (t-t_0)
  \end{array} \right)\xi,
$$
$$
  \eta(2\pi,s,\xi)-\eta(0,s,\xi)=\left(\pi-(\pi/4)(\xi_1^2+\xi_2^2)\right)\left(
  \begin{array}{c}
    \xi_1 \\ \xi_2
  \end{array} \right),
$$
and so (A1) and (A2) of Theorem 2 are satisfied. \qed

\vskip0.4truecm
\noindent
On the basis of the approach presented in Section 2 we now perturb the van der
Pol
equation by means of !a! higher order (with respect to $\varepsilon$)
non-autonomous
term to obtain, for $\varepsilon>0$ sufficiently small,
an integral operator with nonzero topological degree in a suitable open set.

\noindent
In fact, consider the following perturbed system obtained from (\ref{vdp}) by
adding the !forcing! term $\varepsilon\sqrt{\varepsilon}\sin
\left(\frac{t}{1+\varepsilon\sqrt{\varepsilon}\mu}\right)$:
\begin{equation}\label{vdpp}
  \ddot x-\varepsilon(1-x^2)\dot x+x+\varepsilon\sqrt{\varepsilon}\sin
\left(\frac{t}{1+\varepsilon\sqrt{\varepsilon}\mu}\right)=0.
\end{equation}

Consider the following operator $F_{\varepsilon, \mu}: {\rm
C}\left([0,2\pi],{\rm R}^2\right) \to
{\rm C}\left([0,2\pi],{\rm R}^2\right)$  associated to (\ref{vdpp}):
\begin{eqnarray}\label{cov}
  & & F_{\varepsilon, \mu}(x)(t) =
x(t)-x\left(\frac{2\pi}{1+\varepsilon\sqrt{\varepsilon}\mu}\right)-\nonumber\\
  & & -\int\limits_0^t\left(
  \begin{array}{l}
    \qquad \qquad \qquad \quad x_2(\tau)\\
- x_1(\tau)
+\varepsilon(1-x_1^2(\tau))x_2(\tau)-\varepsilon\sqrt{\varepsilon}\sin
\left(\frac{\tau}{1+\varepsilon\sqrt{\varepsilon}\mu}\right)\\   \end{array}
\right) d\tau.
\end{eqnarray}

\noindent
Let $\phi_1(\xi)=\left(
  \begin{array}{l}
    \; \; \quad 0\\
(1-\xi_1^2)\xi_2
\end{array} \right),\ $
$\phi_2(t,\xi,\varepsilon)=\left(
  \begin{array}{l}
   \; \; 0\\
 - \sin t
\end{array} \right),\ $
$T_\mu=\frac{2\pi}{1+\varepsilon\sqrt{\varepsilon}\mu}$
and define the corresponding functions $\eta_1$ and $\eta_2.$

\vskip0.4truecm
\noindent
We are now in the position to formulate the following result.

\vskip0.2truecm
\noindent
{\bf Proposition 2.}
There exists
$\mu_0>0$ such that for every $\mu\in(-\mu_0,\mu_0)$ and $\delta\in(!0!,2)$
there exists
$\varepsilon_0>0$ such that
\begin{eqnarray}\label{p2}
& & {\rm deg}\left(F_{\varepsilon,
\mu},W\left(T_\mu,U_{\max\{0,\delta\}}\right)\backslash \overline
W\left(T_\mu,{U}_{\min\{0,\delta\}}\right) \right) = !1!
\end{eqnarray}
for !any! $\varepsilon\in(0,\varepsilon_0).$

\vskip0.2truecm
\noindent
{\bf Proof.}

\noindent
As already observed
$$
\eta_1(2\pi,s,\xi)-\eta_1(0,s,\xi)=\left(\pi-(\pi/4)(\xi_1^2+\xi_2^2)\right)\left(
  \begin{array}{c}
    \xi_1 \\ \xi_2
  \end{array} \right).
$$
Therefore
\begin{eqnarray}
& & \eta_1(2\pi,s,\xi)-\eta_1(0,s,\xi)=0,\ !\mbox{for any}! \ s\in[0,2\pi] \
!\mbox{and any}! \ \xi\in \partial U_0, \nonumber \\
& & \eta_1(2\pi,s,\xi)-\eta_1(0,s,\xi)\not=0,\ !\mbox{for any}! \ s\in[0,2\pi],\
!\mbox{and any}! \ \xi\in \partial U_{\delta},\nonumber
\end{eqnarray}
moreover
\begin{equation}\label{d1}
  {\rm deg}(\eta_1(2\pi,0,\cdot),U_\delta)=!!1!!
.
\end{equation}

\noindent
It is also easy to see that
$$
  \eta_2(2\pi,s,\xi)-\eta_2(0,s,\xi)=\left(
  \begin{array}{c}
    0 \\ -\frac{3\pi}{4}
  \end{array} \right)
$$
and so
\begin{equation}\label{d2}
  {\rm deg}(\eta_2(2\pi,0,\cdot),U_0)=0.
\end{equation}

\noindent
By (\ref{d1}), (\ref{d2}) and Theorem 4 we obtain (\ref{p2}). \qed

\vskip0.4truecm
\noindent

!So the set $W(T,U_\delta)$ contains the set $X_{\varepsilon,\mu}$ of 
$\frac{2\pi}{1+\varepsilon\sqrt{\varepsilon}\mu}$-periodic solutions 
of the equation (\ref{vdpp}) and the topological degree
of the set $X_{\varepsilon,\mu}$ is not equal to 0 with respect to the integral 
operator $F_{\varepsilon,\mu}$ for the values $\varepsilon,\mu$ 
prescribed above.

From the physical point of view the 
term 
$\varepsilon\sqrt{\varepsilon}\sin
\left(w_{\varepsilon,\mu} t\right)$
means the external voltage including in the van der Pol oscillator 
described by (\ref{vdp}) (see \cite{an}).
So we have proved that 
$w_{\varepsilon,\mu}=\frac{1}{1+\varepsilon\sqrt{\varepsilon}\mu}$
belongs to the frequency pulling range for the forced van der Pol oscillator described by (\ref{vdpp}).
This phenomena is very useful in radio engineering,
but the classical form of the external voltage is $\varepsilon\sin
\left(w_\varepsilon t\right).$
Therefore our approach lets to economize the electricity.!

Observe that, by a suitable choice of the function of the parameter $\mu$ in
(\ref{vdpp}), which
multiplies the term $\varepsilon\sqrt{\varepsilon}$, we can obtain any
prescribed value of the topological
degree in (\ref{p2}).

\end{document}